\documentclass[10pt,draft,twoside,leqno,a4paper]{article}
\usepackage{amsmath,amsthm,amsfonts,amssymb,verbatim}
\newtheorem{Theorem}{Theorem}
\newtheorem{Lemma}{Lemma}

\newcommand{\dash}{\mathchoice
    {\mkern.70mu\raise.50ex\hbox{\vrule height.1ex width.40em depth0pt}}
    {\mkern.40mu\raise.48ex\hbox{\vrule height.1ex width.30em depth0pt}}
    {\mkern.33mu\raise.30ex\hbox{\vrule height.1ex width.25em depth0pt}}
    {\mkern.10mu\raise.20ex\hbox{\vrule height.1ex width.20em depth0pt}}
    }

\newcommand{\R}{\mathbb R}

\renewcommand{\d}{{\mathrm d}}

\newcommand{\Cab}{C(a,b)}

\newcommand{\const}{c}

\renewcommand{\d}{\,{\mathrm d}}

\begin{document}
\title{A sharp weighted Wirtinger inequality}
\author{Tonia Ricciardi\thanks{
Supported in part by 
Regione Campania L.R.~5/02 and by 
the MIUR National Project {\em Variational Methods and
Nonlinear Differential Equations}.}\\
{\small Dipartimento di Matematica e Applicazioni}\\
{\small Universit\`a di Napoli Federico II}\\
{\small Via Cintia, 80126 Naples, Italy}\\
{\small fax: +39 081 675665}\\
{\small\tt{ tonia.ricciardi@unina.it}}\\
}
\date{}
\maketitle
\begin{abstract}
We obtain a sharp estimate for the best constant $C>0$ in the 
Wir\-tin\-ger type inequality
\[
\int_0^{2\pi}\gamma^pw^2\le C\int_0^{2\pi}\gamma^qw'^2
\]
where $\gamma$ is bounded above and below away from zero,
$w$ is $2\pi$-periodic and such that
$\int_0^{2\pi}\gamma^pw=0$, and $p+q\ge0$.
Our result generalizes an inequality of Piccinini and Spagnolo. 
\end{abstract}
%\begin{description}
%\item {\textsc{Key Words:}} linear elliptic equation, measurable coefficients,
%H\"older regularity
%\item {\textsc{MSC 2000 Subject Classification:}} 35J60
%\end{description}
%%%%%%%%%%%%%%%%%%%%%%%%%%%%%%%%%%%%%%%%%%%%%%%%%%%%%%%%%
%%%%%%%%%%%%%%%%%%%%%%%%%%%%%%%%%%%%%%%%%%%%%%%%%%%%%%%%%
%   wirtinger
%%%%%%%%%%%%%%%%%%%%%%%%%%%%%%%%%%%%%%%%%%%%%%%%%%%%%%%%%
%%%%%%%%%%%%%%%%%%%%%%%%%%%%%%%%%%%%%%%%%%%%%%%%%%%%%%%%%
Let $\Cab>0$ denote the best constant in the
following weighted Wirtinger type
inequality:
\begin{align}
\label{wirtinger}
\int_0^{2\pi}a\,w^2\le\Cab\int_0^{2\pi}bw'^2,
\end{align}
where $w\in H_{\mathrm{loc}}^1(\R)$ is $2\pi$-periodic and satisfies
the constraint
\begin{equation}
\label{constraint}
\int_0^{2\pi}aw=0,
\end{equation}
and $a,b\in\mathcal B$ with
\[
\mathcal B=\{a\in L^\infty(\R)\ :\ a\ \mathrm{is}\ 2\pi\mathrm{-periodic\ and\ }\inf a>0\}.
\]
Here and in what follows, for every measurable function $a$
we denote by $\inf a$ and $\sup a$
the essential lower bound and the essential upper bound of $a$,
respectively.
For every $L>1$, we denote
\[
\mathcal B(L)=\left\{a\in L^\infty(0,2\pi)\ :\ a\ \mathrm{is}\ 2\pi\mathrm{-periodic,}\ 
\inf a=1\ \mathrm{and}\
\sup a=L\right\}.
\]
Our aim in this note is to prove:
\begin{Theorem}
\label{thm:wirtinger}
Suppose $a=\gamma^p$ and $b=\gamma^q$ for some $\gamma\in\mathcal B(M)$,
$M>1$, 
and for some $p,q\in\R$ such that $p+q\ge0$.
Then
\begin{equation}
\label{wpq}
C(\gamma^p,\gamma^q)\le\left(
\frac{\frac{1}{2\pi}\int_0^{2\pi}\gamma^{(p-q)/2}}
{\frac{4}{\pi}\arctan\left(M^{-(p+q)/4}\right)}
\right)^2.
\end{equation}
If $p+q>0$, then equality holds in \eqref{wpq} if and only if
$\gamma(\theta)=\bar\gamma_{p,q}(\theta+\varphi)$ for some
$\varphi\in\R$, where
\begin{equation*}
\bar\gamma_{p,q}(\theta)=\begin{cases}
1,&\mathrm{if}\ 0\le\theta<c_{p,q}\frac{\pi}{2},\ \pi\le\theta<\pi+c_{p,q}\frac{\pi}{2}\\
M,&\mathrm{if}\ c_{p,q}\frac{\pi}{2}\le\theta<\pi,\ \pi+c_{p,q}\frac{\pi}{2}\le\theta<2\pi
\end{cases},
\end{equation*}
with  
\[
c_{p,q}=\frac{2}{1+M^{-(p-q)/2}}.
\]
Furthermore, equality holds in \eqref{wirtinger}--\eqref{constraint}
with $a(\theta)=\bar\gamma_{p,q}^p(\theta+\varphi)$ and 
$b(\theta)=\bar\gamma_{p,q}^q(\theta+\varphi)$ if and only if 
$w(\theta)=\bar w_{p,q}(\theta+\varphi)$
where
\begin{align*}
\bar w_{p,q}&(\theta)=\\
=&\begin{cases}
\sin\left[\sqrt\mu\left(c_{p,q}^{-1}\theta-\frac{\pi}{4}\right)\right],\\
\hspace{5cm}\mathrm{if}\ 0\le\theta<c_{p,q}\frac{\pi}{2}\\
M^{-(p+q)/4}\cos\left[\sqrt\mu
\left(\frac{\pi}{2}+c_{p,q}^{-1}M^{(p-q)/2}(\theta-c_{p,q}\frac{\pi}{2})-\frac{3\pi}{4}\right)\right],\\
\hspace{5cm}\mathrm{if}\ c_{p,q}\frac{\pi}{2}\le\theta<\pi\\
-\sin\left[\sqrt\mu\left(
\pi+c_{p,q}^{-1}(\theta-\pi)-\frac{5\pi}{4}\right)\right],\\
\hspace{5cm}\mathrm{if}\ \pi\le\theta<\pi+c_{p,q}\frac{\pi}{2}\\
-M^{-(p+q)/4}\cos\left[\sqrt\mu\left(\frac{3\pi}{2}+c_{p,q}^{-1}M^{(p-q)/2}
(\theta-\pi-c_{p,q}\frac{\pi}{2})-\frac{7\pi}{4}\right)\right],\\
\hspace{5cm}\mathrm{if}\ \pi+c_{p,q}\frac{\pi}{2}\le\theta<2\pi
\end{cases},
\end{align*}
and $\mu=\left((4/\pi)\arctan M^{-(p+q)}\right)^2$.
\par
If $p+q=0$, then
\eqref{wpq} is an equality for any weight function $\gamma$.
Equality is attained in \eqref{wirtinger}--\eqref{constraint}
with $a=\gamma^p$ and $b=\gamma^{-p}$ if and only if
\[
w(\theta)=C\cos\left(\frac{2\pi}{\int_0^{2\pi}\gamma^p}\int_0^\theta\gamma^p+\varphi\right),
\]
for some $C\neq0$ and $\varphi\in\R$.
\end{Theorem}
Note that when $p=q=0$, Theorem~\ref{thm:wirtinger} yields $C(1,1)=1$
according to the classical Wirtinger inequality.
When $p=q\neq0$, the estimate \eqref{wpq} reduces to the estimate obtained 
by Piccinini and Spagnolo in \cite{PS}.
More related results may be found in \cite{Be,CD,DGS}
and in the references therein.
We begin by recalling in the following lemma
the Wirtinger inequality of Piccinini and Spagnolo~\cite{PS}. 
\begin{Lemma}[\cite{PS}]
\label{lem:PS}
Suppose $b=a\in\mathcal B(L)$. Then,
\begin{equation}
\label{wPS}
C(a,a)\le
\left(\frac{4}{\pi}\arctan L^{-1/2}
\right)^{-2}.
\end{equation}
Equality holds in \eqref{wPS} if and only if 
$a(\theta)=\bar a(\theta+\varphi)$ for some $\varphi\in\R$, 
where $\bar a$ is defined by
\begin{equation}
\label{a1def}
\bar a(\theta)=\begin{cases}
1,&\mathrm{if}\ 0\le\theta<\frac{\pi}{2},\ \pi\le\theta<\frac{3\pi}{2}\\
L,&\mathrm{if}\ \frac{\pi}{2}\le\theta<\pi,\ \frac{3\pi}{2}\le\theta<2\pi
\end{cases}
\end{equation}
and
equality holds in \eqref{wirtinger}--\eqref{constraint}
with $a(\theta)=b(\theta)=\bar a(\theta+\varphi)$ 
if and only if $w(\theta)=\bar w(\theta+\varphi)$, where
\begin{equation}
\bar w(\theta)=\begin{cases}
\sin\left[\sqrt\lambda\left(\theta-\frac{\pi}{4}\right)\right],
&\mathrm{if}\ 0\le\theta<\frac{\pi}{2}\\
L^{-1/2}\cos\left[\sqrt\lambda\left(\theta-\frac{3\pi}{4}\right)\right],
&\mathrm{if}\ \frac{\pi}{2}\le\theta<\pi\\
-\sin\left[\sqrt\lambda\left(\theta-\frac{5\pi}{4}\right)\right],
&\mathrm{if}\ \pi\le\theta<\frac{3\pi}{2}\\
-L^{-1/2}\cos\left[\sqrt\lambda\left(\theta-\frac{7\pi}{4}\right)\right],
&\mathrm{if}\ \frac{3\pi}{2}\le\theta<2\pi
\end{cases},
\end{equation}
where $\lambda=\left(4\pi^{-1}\arctan L^{-1/2}\right)^2$.
\end{Lemma}
In order to prove Theorem~\ref{thm:wirtinger},
we need the following lemma, which yields an estimate for $\Cab$ for 
\textit{arbitrary} weight functions $a,b$.
\begin{Lemma}
\label{lem:prelimCab}
Let $a,b\in\mathcal B$.
The following estimate holds:
\begin{equation}
\label{prelimCab}
\Cab\le
\left(
\frac{\frac{1}{2\pi}\int_0^{2\pi}\sqrt{ab^{-1}}}
{\frac{4}{\pi}\arctan\left({\frac{\inf{ab}}
{\sup{ab}}}\right)^{1/4}}\right)^2.
\end{equation}
If $\sqrt{ab}\in\mathcal B(L)$, $L>1$, then
\begin{equation}
\label{a'b'}
\frac{\Cab}{\left(\frac{1}{2\pi}\int_0^{2\pi}\sqrt{ab^{-1}}\right)^2}
=\sup_{\sqrt{a'b'}\in\mathcal B(L)}
\frac{C(a',b')}{\left(\frac{1}{2\pi}\int_0^{2\pi}\sqrt{a'b'^{-1}}\right)^2}
=\left(\frac{4}{\pi}\arctan L^{-1/2}\right)^{-2}
\end{equation}
if and only if the following equation is satisfied: 
\begin{equation}
\label{functionaleq}
a(\theta(\tau))b(\theta(\tau))=\bar a^2(\tau+\varphi)
\qquad\mathrm{a.e.}\ \tau\in(0,2\pi),
\ \mathrm{for\ some\ }\varphi\in\R,
\end{equation}
where $\theta(\tau)$ is the homeomorphism of $\R$
defined by
\begin{equation}
\label{taudef}
\tau(\theta)=\frac{1}{\const}\int_0^\theta
\sqrt{\frac{a(\tilde\theta)}{b(\tilde\theta)}}\d\tilde\theta,
\end{equation}
$\const$ is defined by
\begin{equation}
\label{const}
\const=\frac{1}{2\pi}\int_0^{2\pi}
\sqrt{\frac{a(\tilde\theta)}{b(\tilde\theta)}}\d\tilde\theta,
\end{equation}
and $\bar a$ is the function defined in Lemma~\ref{lem:PS}.
\par
If $b=a^{-1}$, then $C(a,a^{-1})=\left((2\pi)^{-1}\int_0^{2\pi}a\right)^2$
and equality is attained in \eqref{wirtinger}--\eqref{constraint} with $b=a^{-1}$
if and only if 
$w(\theta)=C\cos(2\pi(\int_0^{2\pi}a)^{-1}\int_0^\theta a+\varphi)$
for some $C\neq0$ and $\varphi\in\R$.
\end{Lemma}
\begin{proof}
Under the change of variables $\theta=\theta(\tau)$ 
defined by \eqref{taudef}--\eqref{const},
setting $\alpha(\tau)=a(\theta(\tau))$,
$\beta(\tau)=b(\theta(\tau))$, $\xi(\tau)=w(\theta(\tau))$,
we obtain
\begin{align*}
&\alpha\theta'=\const\sqrt{\alpha\beta},
&&\beta{\theta'}^{-1}=\const^{-1}\sqrt{\alpha\beta},
\end{align*}
and therefore:
\begin{align*}
\int_0^{2\pi}aw^2\d\theta=&\int_0^{2\pi}\alpha\theta'\xi^2\d\tau
=\const\int\sqrt{\alpha\beta}\xi^2\d\tau\\
\int_0^{2\pi}aw\d\theta=&\int_0^{2\pi}\alpha\theta'\xi\d\tau
=\const\int\sqrt{\alpha\beta}\xi\d\tau=0\\
\int_0^{2\pi}bw'^2\d\theta=&\int_0^{2\pi}\beta{\theta'}^{-1}
\xi'^2\d\tau
=\const^{-1}\int\sqrt{\alpha\beta}\xi'^2\d\tau.
\end{align*}
Upon substitution, 
\eqref{wirtinger}--\eqref{constraint} takes the form:
\begin{equation}
\label{subswirtinger}
\int_0^{2\pi}\sqrt{\alpha\beta}\xi^2\d\tau
\le\frac{\Cab}{\left(\frac{1}{2\pi}\int_0^{2\pi}\sqrt{ab^{-1}}\right)^2}
\int_0^{2\pi}\sqrt{\alpha\beta}\xi'^2\d\tau,
\end{equation}
with constraint
\begin{equation}
\label{subsconstraint}
\int_0^{2\pi}\sqrt{\alpha\beta}\xi\d\tau=0.
\end{equation}
If $\sqrt{ab}\in\mathcal B(L)$, in view of Lemma~\ref{lem:PS} we obtain
\begin{align}
\label{Calphabeta}
\frac{\Cab}{\left(\frac{1}{2\pi}
\int_0^{2\pi}\sqrt{ab^{-1}}\right)^2}
=C(\sqrt{\alpha\beta},\sqrt{\alpha\beta})
\le&\left(\frac{4}{\pi}\arctan
\sqrt{\frac{\inf\sqrt{\alpha\beta}}{\sup\sqrt{\alpha\beta}}}\right)^{-2}\\
\nonumber
=&\left(\frac{4}{\pi}\arctan
\left({\frac{\inf{ab}}{\sup{ab}}}\right)^{1/4}\right)^{-2}.
\end{align}
This yields \eqref{prelimCab}. 
Moreover, we have
$C(\sqrt{\alpha\beta},\sqrt{\alpha\beta})
=\left((4/\pi)\arctan L^{-1/2}\right)^{-2}$
if and only if
$\sqrt{\alpha(\tau)\beta(\tau)}=\bar a(\tau+\varphi)$,
for some $\varphi\in\R$. 
That is, 
\eqref{a'b'} holds if and only if \eqref{functionaleq}
holds.
\par
If $b=a^{-1}$, then \eqref{subswirtinger}--\eqref{subsconstraint}
takes the form
\[
\int_0^{2\pi}\xi^2\d\tau\le\frac{C(a,a^{-1})}
{\left(\frac{1}{2\pi}\int_0^{2\pi}a\right)^2}
\int_0^{2\pi}\xi'^2\d\tau
\]
with constraint
\[
\int_0^{2\pi}\xi\d\tau=0.
\]
Therefore, by the classical Wirtinger inequality,
\[
C(a,a^{-1})=\left(\frac{1}{2\pi}\int_0^{2\pi}a\right)^2
\]
and equality holds in \eqref{wirtinger}--\eqref{constraint}
with $b=a^{-1}$ if and only if $\xi(\tau)=C\cos(\tau+\varphi)$
for some $C\neq0$ and $\varphi\in\R$, that is, if and only if
$w(\theta)=C\cos(2\pi(\int_0^{2\pi}a)^{-1}\int_0^\theta a+\varphi)$,
as asserted.
\end{proof}
\begin{Lemma}
\label{lem:functionaleq}
Suppose $a,b$ satisfy $\sqrt{ab}\in\mathcal B(L)$,
$L>1$,
and \eqref{functionaleq},
where $\theta(\tau)$ is defined in \eqref{taudef} and $c$ is defined by \eqref{const}.
Suppose
\begin{equation}
\label{apcond}
a=\gamma^p,\qquad b=\gamma^q
\end{equation}
for some $\gamma\in\mathcal B(M)$, with $M=L^{2/(p+q)}$, 
and for some $p,q\in\R$ such that $p+q>0$.
Then $\gamma(\theta)=\bar\gamma_{p,q}(\theta+\varphi)$
for some $\varphi\in\R$, where $\bar\gamma_{p,q}$ is the function
defined in Theorem~\ref{thm:wirtinger}.
\end{Lemma}
\begin{proof}
When $p+q>0$, we have $\gamma^{(p+q)/2}\in\mathcal B(L)$.
In view of \eqref{functionaleq} and \eqref{apcond} we have
\[
\gamma(\theta(\tau))=\bar a^{2/(p+q)}(\tau+\psi),
\qquad\forall\tau\in\R
\]
for some $\psi\in\R$.
It follows that
\begin{align}
\label{thetatau}
\theta(\tau)
=c\int_0^{\tau}\sqrt{\frac{b(\theta(\bar\tau))}{a(\theta(\bar\tau))}}\d\bar\tau
=c\int_0^\tau\bar a^{-(p-q)/(p+q)}(\bar\tau+\psi)\d\bar\tau
\end{align}
and, in view of the $2\pi$-periodicity of $a$ and $b$,
\[
c=\left(\frac{1}{2\pi}\int_0^{2\pi}
\sqrt{\frac{b(\theta(\bar\tau))}{a(\theta(\bar\tau))}}\d\bar\tau\right)^{-1}
=\left(\frac{1}{2\pi}\int_0^{2\pi}\bar a^{-(p-q)/(p+q)}(\bar\tau)\d\bar\tau\right)^{-1}.
\]
Setting 
\[
h_{p,q}(\tau)=c\int_0^\tau\bar a^{-(p-q)/(p+q)}(\bar\tau)\d\bar\tau,
\]
we have $\theta(\tau-\psi)=h_{p,q}(\tau)-h_{p,q}(\psi)$
for every $\tau\in\R$,
and consequently $\tau(\theta)=h_{p,q}^{-1}(\theta+h_{p,q}(\psi))-\psi$.
In view of the definition of $\bar a$ with $L=M^{(p+q)/2}$, 
we have:
\begin{align*}
\int_0^\tau\bar a^{-(p-q)/(p+q)}&(\bar\tau)\d\bar\tau=\\
=&\begin{cases}
\tau,&\mathrm{if}\ 0\le\tau<\frac{\pi}{2}\\
\frac{\pi}{2}+M^{-(p-q)/2}(\tau-\frac{\pi}{2}),
&\mathrm{if}\ \frac{\pi}{2}\le\tau<\pi\\
\frac{\pi}{2}(1+M^{-(p-q)/2})+\tau-\pi,
&\mathrm{if}\ \pi\le\tau<\frac{3\pi}{2}\\
\frac{\pi}{2}(2+M^{-(p-q)/2})+M^{-(p-q)/2}(\tau-\frac{3\pi}{2}),
&\mathrm{if}\ \frac{3\pi}{2}\le\tau<2\pi
\end{cases}.
\end{align*}
In particular,
we derive
\[
c=\frac{2}{1+M^{-(p-q)/2}}=c_{p,q}.
\]
It follows that $h_{p,q}(\tau)$ is the piecewise linear homeomorphism
of $\R$ defined in $[0,2\pi)$ by
\[
h_{p,q}(\tau)=\begin{cases}
c_{p,q}\tau,
&\mathrm{if}\ 0\le\tau<\frac{\pi}{2}\\
c_{p,q}\left[\frac{\pi}{2}+M^{-(p-q)/2}(\tau-\frac{\pi}{2})\right],
&\mathrm{if}\ \frac{\pi}{2}\le\tau<\pi\\
c_{p,q}\left[\frac{\pi}{2}(1+M^{-(p-q)/2})+\tau-\pi\right],
&\mathrm{if}\ \pi\le\tau<\frac{3\pi}{2}\\
c_{p,q}\left[\frac{\pi}{2}(2+M^{-(p-q)/2})+M^{-(p-q)/2}(\tau-\frac{3\pi}{2})\right],
&\mathrm{if}\ \frac{3\pi}{2}\le\tau<2\pi
\end{cases}
\]
and by $h_{p,q}(\tau+2\pi n)=2\pi n+h_{p,q}(\tau)$,
for any $\tau\in[0,2\pi)$ and for any integer $n$.
Inversion yields
\[
h_{p,q}^{-1}(\theta)=\begin{cases}
c_{p,q}^{-1}\theta,
&\mathrm{if}\ 0\le\theta<c_{p,q}\frac{\pi}{2}\\
\frac{\pi}{2}+c_{p,q}^{-1}M^{(p-q)/2}(\theta-c_{p,q}\frac{\pi}{2}),
&\mathrm{if}\ c_{p,q}\frac{\pi}{2}\le\theta<\pi\\
\pi+c_{p,q}^{-1}(\theta-\pi),
&\mathrm{if}\ \pi\le\theta<\pi+c_{p,q}\frac{\pi}{2}\\
\frac{3\pi}{2}+c_{p,q}^{-1}M^{(p-q)/2}(\theta-\pi-c_{p,q}\frac{\pi}{2}),
&\mathrm{if}\ \pi+c_{p,q}\frac{\pi}{2}\le\theta<2\pi
\end{cases},
\]
for $\theta\in[0,2\pi)$ and $h_{p,q}^{-1}(\theta+2\pi n)=2\pi n+h_{p,q}^{-1}(\theta)$
for any $\tau\in[0,2\pi)$ and for any integer $n$.
Substitution yields
$\gamma(\theta)=\bar a^{2/(p+q)}\left(h_{p,q}^{-1}(\theta+h_{p,q}(\psi))\right)
=\bar a^{2/(p+q)}\left(h_{p,q}^{-1}(\theta+\varphi)\right)
=\bar\gamma_{p,q}(\theta+\varphi)$, with $\varphi=h_{p,q}(\psi)$.
\end{proof}
Now we can prove Theorem~\ref{thm:wirtinger}.
\begin{proof}
[Proof of Theorem~\ref{thm:wirtinger}]
Estimate~\eqref{prelimCab} with $a=\gamma^p$
and $b=\gamma^q$ yields \eqref{wpq}.
Suppose $p+q>0$. 
In view of Lemma~\ref{lem:prelimCab} and Lemma~\ref{lem:functionaleq}
we have
\[
\frac{C(\gamma^p,\gamma^q)}
{\left(\frac{1}{2\pi}\int_0^{2\pi}\gamma^{(p-q)/2}\right)^2}
=\left(\frac{4}{\pi}\arctan M^{-(p+q)/4}\right)^{-2}
\]
if and only if $\gamma(\theta)=\bar\gamma_{p,q}(\theta+\varphi)$
for some $\varphi\in\R$.
Equality is attained in \eqref{wirtinger}--\eqref{constraint}
with $a(\theta)=\bar\gamma_{p,q}^p(\theta+\varphi)$
and $b(\theta)=\bar\gamma_{p,q}^q(\theta+\varphi)$ if and only if
$w(\theta)=\bar w_{p,q}(\theta+\varphi)$.
\par
If $p+q=0$, then the conclusion follows by Lemma~\ref{lem:prelimCab}
with $a=\gamma^p$ and $b=\gamma^{-p}$.
\end{proof}
\section*{Acknowledgements}
I am grateful to Professor Carlo Sbordone for
many useful and stimulating discussions.
%%%%%%%%%%%%%%%%%%%%%%%%%%%%%%%%%%%%%%%%%%%%%%%%%%%%%%%%%
%%%%%%%%%%%%%%%%%%%%%%%%%%%%%%%%%%%%%%%%%%%%%%%%%%%%%%%%%
%   References
%%%%%%%%%%%%%%%%%%%%%%%%%%%%%%%%%%%%%%%%%%%%%%%%%%%%%%%%%
%%%%%%%%%%%%%%%%%%%%%%%%%%%%%%%%%%%%%%%%%%%%%%%%%%%%%%%%%


\begin{thebibliography}{99}
\bibitem{Be}
P.R.~Beesack,
Integral inequalities of the Wirtinger type, 
Duke Math.\ Jour.\ \textbf{25} (1958), 477--498.
\bibitem{CD}
G.~Croce and B.~Dacorogna,
On a generalized Wirtinger inequality, 
Discrete Cont.\ Dynam.\ Systems \textbf{9} No.~5 (2003), 1329--1341.
\bibitem{DGS}
B.~Dacorogna, W.~Gangbo and N.~Sub\'ia,
Sur une g\'en\'eralisation de l'in\'egalit\'e de Wirtinger, 
Ann.\ Inst.\ H.~Poincar\'e Anal.\ Non Lin\'eaire \textbf{9} (1992), 29--50.
%\bibitem{IS}
%T.~Iwaniec and C.~Sbordone,
%Quasiharmonic fields, 
%Ann.\ Inst.\ H.~Poincar\'e Anal.\ Non Lin\'eaire \textbf{18} No.~5 (2001), 519--572.
\bibitem{PS} 
L.C.~Piccinini and S.~Spagnolo, 
On the H\"older continuity of solutions of second order elliptic equations
in two variables, Ann.\ Scuola Norm.\ Sup.\ Pisa \textbf{26} No.~2 (1972), 391--402.
%\bibitem{R11} 
%T.~Ricciardi,  
%A sharp H\"older estimate for elliptic equations
%in two variables, 
%Proc.\ Roy.\ Soc.\ Edinburgh A, to appear.
%Preprint on arXiv:math.AP/0406495.
\end{thebibliography}
\end{document}